\documentclass[12pt]{article}
\usepackage{latexsym}
\usepackage{amssymb}
\begin{document}
\newtheorem{theorem}{Theorem}[section]
\newtheorem{lemma}[theorem]{Lemma}
\newtheorem{proposition}[theorem]{Proposition}
\newtheorem{corollary}[theorem]{Corollary}
\newtheorem{definition}{Definition}
\newtheorem{question}{Question}
\newtheorem{conjecture}{Conjecture}
\newcommand{\F}{\ensuremath{\mathbb F}}
\newcommand{\N}{\mathcal N}
\newcommand{\R}{\mathcal R}
\newcommand{\Z}{\mathbb Z}
	
\title{Conjugacy classes of maximal cyclic subgroups and nilpotence class of $p$-groups}
\author{M. Bianchi, R.D. Camina, \& Mark L. Lewis }
\date \
\maketitle
	
\begin{abstract} 
In this paper, we set $\eta (G)$ to be the number of conjugacy classes of maximal cyclic subgroups of $G$.  We prove that if $G$ is a $p$-group of order $p^n$ and nilpotence class $l$, then $\eta (G)$ is bounded below by a linear function in $n/l$. \\[1ex]
{\it Keywords:} group covering, nilpotence class, group exponent\\[1ex]
{\it 2020 Mathematics Subject Classification:} 20D15
\end{abstract}
	
\maketitle
	
\section{Introduction}
	
Unless otherwise stated, all groups in this paper are finite.  Motivated by the research area of classifying spaces of families of subgroups of infinite groups, von Puttkamer asks \cite[Question 5.0.9]{von}: does the number of conjugacy classes of maximal cyclic subgroups of a finite $p$-group for a prime $p > 2$ grow with the order of the group?   Clearly, the question is considering noncyclic groups.   More precisely, one can ask whether the number of conjugacy classes of maximal cyclic subgroups of a noncyclic $p$-group of order $p^n$ is at least $n$ when $p$ is an odd prime.  This work was initiated by X. Wu asking the second author this question in a private communication.
	
The question can be phrased in terms of the size of a particular type of covering of a group. Recall, a {\it covering } of a group $G$ is a set of proper subgroups $\{ H_i\}$, called {\it components}, such that $G \subseteq \bigcup_i H_i$.  A covering is called {\it irredundant} if removing any component means the set is no longer a covering.  A {\it normal covering} of $G$ is a covering that is invariant under conjugation by $G$.  Note that the components of a normal covering are not normal, but refers to the fact that the components form orbits under conjugacy.  Often, one takes only one representative of each conjugacy class of components.  Our question above considers normal coverings where all the components of the cover are cyclic.  It is not difficult to see that the set of maximal cyclic subgroups is the only irredundant covering by cyclic subgroups.
	
With this in mind, we introduce the following definition.  Let $G$ be a finite group. Denote the number of conjugacy classes of maximal cyclic subgroups of $G$ by $\eta(G)$.  We note that the questions of von Puttkamer and Wu have negative answers.  In \cite{pre}, the second and third authors along with Yiftach Barnea and Mikhail Ershov find for every prime $p > 3$ and every integer $n \ge 3$, groups of order $p^n$ that all have $\eta = p+2$.  For $p = 2$, we also find that $\eta (G) = 3$ when $G$ is a dihedral $2$-group, a generalized quaternion $2$-group, or a semi-dihedral $2$-group.   For $p = 3$, there is an infinite family of $3$-groups with $\eta = 9$.  See Theorem 7.12 of \cite{pre}.
	
On the other hand, we believe it is extremely rare for this to occur.  In particular, we prove that if we fix the nilpotence class of a $p$-group $G$, then $\eta (G)$ will grow proportionally to $\log_p (|G|)$.  The theorem we prove is the following.
	
\begin{theorem} \label {nilp cl}
Let $G$ be a noncyclic $p$-group of nilpotence class $l \ge 1$ and order $p^n$.  Then $\eta (G) \ge (p-1)(n/l - 2) + p + 1$.
\end{theorem}  

The authors would like to thank Emanuele Pacifici for a number of helpful conversations while working on this paper.  The first author is partially supported by INDAM-GNSAGA.  We would like to thank Avinoam Mann for pointing out \cite{east}.

\section{Abelian groups}

As noted in the introduction, when considering abelian groups calculating $\eta(G)$ is the same as calculating the maximum number of subgroups in an irredundant cyclic covering of $G$. In \cite[Proposition 6 (ii)]{rogerio} the author gives a formula, which involves the Euler totient function, for calculating this number. 

We take a slightly different approach; so we are including our results.  We will show for abelian $p$-groups that a function in $\eta (G)$ gives a lower bound for $|G|$.   We begin by computing $\eta (G)$ when $G$ is the direct product of two cyclic $p$-groups. 

We consider the set $G^{\{p\}} = \{ g^p \mid g \in G \}$.  For an element $g \in G$, it is not difficult to see that that $\langle g \rangle$ is a maximal cyclic subgroup of $G$ if and only if $g \in G \setminus G^{\{p\}}$.

\begin{lemma}\label{two}
Suppose $G \cong C_{p^a} \times C_{p^b}$ with $a \geq b$. Then
	$$ \eta(G) =  p^{(b-1)}((a-b) (p-1) + p +  1) \geq a+b.$$
\end{lemma}

{\bf Proof.}  Let $G = \langle x \rangle \times \langle y \rangle$ with $x$ of order $p^a$ and $y$ of order $p^b$. 
Let $C$ be a cyclic subgroup of $G$.  Either $C$ is a
subgroup of  $\langle (1,y) \rangle$ or there exist integers $n$ and $c$ so that $C$ is
generated by $(x^{p^n},y^c)$ for $0 \leq n < a$ and $0 \leq c < p^b$.
We claim that $C$ is maximal exactly when $n = 0$, or $c$ is relatively
prime to $p$, or $C = \langle (1,y) \rangle$.  To see this, observe that in each of these cases, a generator of $C$ does not lie in $G^{\{p\}}$, so $C$ is maximal.  
Furthermore, we claim that the generator for a maximal
cyclic subgroup is unique if we make the additional restriction that $0 \leq c < \min(p^b, p^{a-n})$.
First note for an integer
$l$ that  $(x^{p^n}, y^c)^{1 + lp^{a-n}} = (x^{p^n}, y^{c + clp^{a-n}})$. 
Since $c$ is coprime to $p$, as $l$ runs through the integers modulo $p^b$, then also $cl$ 
will run through all of ${\mathbb Z}_{p^b}$.  So we are getting as the exponents for
$y$, all elements in the coset $c + p^{a-n} {\mathbb Z}_{p^b}$ in ${\mathbb Z}_{p^b}$.  It is not
difficult to see that $\{0,1,...,p^{a-n} -1\}$ is a transversal for $p^{a-n}{\mathbb Z}_{p^b}$
in ${\mathbb Z}_{p^b}$.  Hence, there exists an integer $c'$ with $0 \le c' \le p^{a-n}-1$
so that $c' + p^{a-n}{\mathbb Z}_{p^b} = c + p^{a-n}{\mathbb Z}_{p^b}$.  It follows that
$(x^{p^n},y^{c'})$ will lie in $\langle (x^{p^n},y^c) \rangle$ and since it has the same
order, it will be a generator.

So, we claim
the following are distinct maximal cyclic subgroups of $G$,\\
$$\langle (x, y^c) \rangle\;\; {\rm for} \;\;0 \leq c \leq p^{b}-1\;\; {\rm and}\;\; \langle (1, y) \rangle$$
$$\langle (x^{p^n}, y^c) \rangle\;\;\;{\rm for}\;\;\; 1 \leq n < a,\;\;\; 0 \leq c <  \min(p^b, p^{a-n})\;\;\; {\rm and}\;\; {\rm gcd}(c,p) =1.$$
To count the number of subgroups of the form $\langle (x^{p^n}, y^c) \rangle$ we consider two cases, when $1 \leq n \leq a-b$
and when $a-b+1 \leq n < a$. For $1 \leq n \leq a-b$ we have $\phi(p^b) = p^{b-1}(p-1)$ such subgroups where $\phi$
is the Euler totient function. For $a-b+1 \leq n \leq a-1$ 
we have $\phi (p^{a-n})$ maximal
cyclic subgroups of the form $ \langle x^{p^n},y^c \rangle$ where $p$ does not divide $c$ and
$1 \leq c \leq p^{a-n}$. With this range on $b$, we obtain $1 \leq a-n \leq b-1$.
Observe that $\phi (p^{a-n}) = p^{a-n-1} (p-1)$.  We can view this as
$p^{i-1} (p-1)$ for $i$ running from 1 to $b-1$ or $p^i (p-1)$ for $i$ running
from 0 to $b-2$.
So, in total we have 
$$p^b + (a-b)p^{(b-1)} (p-1)  +  \sum_{i=0}^{b-2} p^i
(p-1) + 1$$
$$   =  p^b + (a-b) p^{(b-1)} (p-1) + [(p^{(b-1)} - 1)/(p-1)] (p-1) + 1 $$
$$ =  p^b + (a-b) p^{(b-1)} (p-1) + p^{(b-1)}$$
$$ =  p^{(b-1)} ((a - b)(p-1) + p + 1)$$
maximal cyclic subgroups, as desired. A simple inductive argument shows that $p^{(b-1)}((a-b)(p-1) + p + 1) \geq a+b$.
$\Box$\\

It is useful to have a function which encodes this value.  Let $p$ be a prime, and let $a$ and $b$ be positive integers.  We take $k = {\rm max} (a,b)$ and $l = {\rm min} (a,b)$.  We set $g_p (a,b) = p^{(l-1)}((k-l) (p-1) + p +  1)$, so $g_p (a,b) = \eta (G)$ when $G = C_{p^a} \times C_{p^b}$.  In this next lemma, we obtain the lower bound that we need to prove the general lower bound for abelian groups.

\begin{lemma} \label{compute g_p}
Let $a$ and $b$ be positive integers and $p$ a prime.  Then $g_p (a,b) \ge (p-1)(a+b-2) + p + 1 \geq a + b + 1$.
\end{lemma}

{\bf Proof.}
Without loss of generality, we may assume $a \ge b$.  Set $n = a+b$, so that $a = n - b$.  We work to show $g_p (a,b) \ge (p-1)(n-2) + p + 1 \geq n + 1$.  We have $1 \le b \le n/2$.  We see that $g_p (n-b,b) = p^{b-1} ((n-2b)(p-1) + p+1)$.  We can view $n$ as fixed, and this becomes a function in one variable:
$$f(x) = g_p (n-x,x) = p^{x-1}((n-2x)(p-1) + p+1).$$  We want to find the minimal value for $f(x)$ on $[1,n/2]$.  We use calculus to see that
$$f'(x) = p^{x-1}(((\ln p)(n-2x) - 2)(p-1) + (\ln p)(p+1)).$$  
Setting this equal to $0$, we obtain the critical value of
$$x = n/2 + (p+1)/(2(p-1)) - 1/(\ln p).$$ 
We claim that $(p+1)/(2(p-1)) - 1/(\ln p)$ is positive for all primes $p$.  (This can be shown using calculus or graphing using a computer.)  Hence, this critical point is not in our interval.

Now, $f(1) = ((n-2)(p-1)) + p +1 = p(n-1) -n +3$ and $f(n/2) = p^{n/2 - 1}(p+1)$.  Recall that Fermat's theorem tells us that any local extreme values would occur when the derivative was zero.  Since the derivative of this function is never zero on this interval, we know that our minimum is the smallest of these two values. Set $f_1(x) = ((x-2)(p-1)) + p + 1$ and $f_{n/2}(x) = p^{x/2 -1}(p+1)$ for $x \geq 2$.  Then $f_1(2) = f_{n/2}(2)$, but $f'_{n/2} (x) > f'_{1} (x)$; so $f_{n/2}(x) > f_{1}(x)$ for $x> 2$. Thus, $g_p (n-b,b) \geq (p-1)(n-2) + p + 1$.  Finally, note that $(p-1)(n-2) + p + 1 \geq n-2 + 3= n+1.$
$\Box$\\

We also need a refinement on the value of $\eta$ of a direct product of an abelian $p$-group with a cyclic $p$-group. 

\begin{lemma}\label{direct} 
Suppose $G \cong H \times C_{p^a}$ where $H$ is an abelian $p$-group, then $\eta(G) \geq (a+1)\eta(H) + 1$.
\end{lemma}

{\bf Proof.}
We write $Y = \langle y \rangle \cong C_{p^a}$, so that $G = H \times Y$.  Let $\{ h_i \}$, for $1 \leq i \leq \eta (H)$, be representatives of generators for the maximal cyclic subgroups of $H$.  We now prove that $\langle (h_i, y^{p^c})\rangle$ for $0 \le c \le a$ and $1 \leq i \leq \eta(H)$, along with $\langle (1,y)\rangle$ are distinct maximal cyclic subgroups of $G$. This gives the required count for the lower bound.

First note that $h_i \in H \setminus H^{\{p\} }$, since $\langle h_i \rangle$ is maximal in $H$, and thus, $(h_i, 1)$, $(h_i, y^{p^c})$, and $(1, y)$ all lie in $G \setminus G^{\{p\} }$ for $c \in \{ 0, \dots, a \}$.  (It is not difficult to see that $G^{\{p\} } = \{ (\alpha,\beta) \mid \alpha \in H^{\{ p\} }, \beta \in \langle y^p \rangle \}$.)  Thus, each of these elements generate a maximal cyclic subgroup of $G$.  We need to show that they generate different cyclic subgroups of $G$. Suppose $\langle(h_i, y^{p^c})\rangle = \langle(h_j, y^{p^d})\rangle$ for $0 \le c,d \le a$.  In both cases, by projecting into $H$, it is not difficult to see that $\langle h_i\rangle = \langle h_j\rangle$, and thus, $i = j$.  We are left to consider $\langle(h_i, y^{p^c})\rangle = \langle (h_i, y^{p^d})\rangle$.   Projecting into $Y$, this forces $\langle y^{p^c}\rangle = \langle y^{p^d} \rangle$, and thus, $c=d$ as required. 
%
$\Box$ \\

We now obtain our lower bound on $\log_p (|G|)$ when $G$ is an abelian $p$-group in terms of $\eta (G)$.

\begin{theorem}\label{one} 
Let $G$ be an abelian, noncyclic group of order $p^n$.  Then $$\eta (G) \ge (p-1)(n-2) + p + 1 \geq n + 1.$$
\end{theorem}

{\bf Proof.} We work by induction on $|G|$.  Suppose first that $G = C_{p^a} \times C_{p^b}$ for positive integers $a$ and $b$.  By Lemma \ref{two}, we have $\eta (G) = g_p (a,b)$ and applying Lemma \ref{compute g_p}, we have $g_p (a,b) \ge (p-1)(a+b-2) + p + 1 \ge a+b + 1$.

Now suppose $H$ is not cyclic and $|H| = p^h$ and $G = H \times C_{p^a}.$  We have by induction, $\eta (H) \ge (p-1)(h-2) + p + 1$.

By Lemma \ref{direct}, 
\begin{eqnarray*}
	\eta (G) = \eta (H \times C_{p^a}) & \ge & (a + 1) \eta (H)+ 1 \\ 
	& \ge & (a + 1)((h-2)(p-1) + p + 1) + 1 \\
	& = & (a+1)(h-2)(p-1) + a (p+1) + p + 1 + 1 \\
	& = & (a+1)(h-2)(p-1) + a (p-1) + 2a + p + 2 \\
	& = & ((a+1)(h-2) + a)(p-1) + 2a + p + 2  \\
	& \ge &  (h - 2 + a)(p-1) + p + 1 
\end{eqnarray*} since $n = h + a$, this gives the result.

Finally, note that $(p-1)(n-2) + p + 1 \geq n-2 + 3= n+1.$ $\Box$\\[1ex]

When $p$ is odd, we can improve the inequality in Theorem \ref{one}.

\begin{corollary}
If $p$ odd and $G$ is a noncyclic, abelian group of order $p^n$, then $\displaystyle \eta(G) \geq (\frac{p+1}{2})n$.
\end{corollary}

{\bf Proof:}  We apply the result from Theorem \ref{one}:
\begin{eqnarray*}
	(p-1)(n-2) + p + 1 & = & (p-1)n -2p + 2 + p + 1 \\
	& = & (p+1)n/2 + (p-3)n/2 - p + 3\\
	& = & (p+1)n/2 + (p-3)(n-2)/2.
\end{eqnarray*}  Since $p \ge 3$ and $n \ge 2$, this is at least
$(p+1)n/2$. 
$\Box$\\





\section{$\eta$ and nilpotence class}

We now show that if we fix the nilpotence class of a $p$-group $G$, then a function of $\eta (G)$ gives a lower bound for $|G|$.  We will see that this can be viewed as generalization of Theorem \ref{one} for abelian groups.  We define the terms of the lower central series of $G$ inductively as follows, $G^1 = G$ and for $i \ge 1$, we set $G^{i+1} = [G^{i},G]$.  Observe that $G$ is nilpotent of nilpotence class $l$ if $G^{l+1} = 1$ and $G^l > 1$.  

We first prove a lemma regarding the order of elements of $G^l$ when $G$ has nilpotence class $l$.  Taking $|G| = p^n$, note that $n \ge l+1$, so $\lfloor n/(l+1) \rfloor \ge 1$.  It has been suggested that this next lemma is related to \cite{east}.

\begin{lemma} \label{orders}
Let $G$ be a $p$-group of nilpotence class $l \ge 2$ and order $p^n$.  Then $G^l$ has exponent dividing $p^{\lfloor n/(l+1) \rfloor}$.
\end{lemma}

{\bf Proof.}
We work by induction on $l$.  We begin with the case that $l = 2$.  If $|G^2| \le p^{n/3}$, then we have the result.  Thus, we may assume that $|G^2| > p^{ n/3}$ and so, $|G/G^2| < p^{2n/3}$.  We can find $a_1, \dots, a_k \in G$ so that $G/G^2 = \langle a_1 G^2 \rangle \times \langle a_2 G^2 \rangle \times \cdots \times \langle a_k G^2 \rangle$.  We know that $G^2$ is central, and so it is abelian.  Also, it is generated by $\{ [a_i,a_j] \mid 1 \le i < j \le k \}$.  Since $G$ is a $p$-group, it suffices to show that $o([a_i,a_j])$ is less than or equal to $p^{n/3}$ for all $1 \le i < j \le k$.   Observe that $o(a_iG^2)o(a_jG^2)$ divides $|G/G^2| < p^{2n/3}$, and so, without loss of generality $o(a_iG^2)$ is less than or equal to $p^{n/3}$.  Now, $[a_i,a_j]^{\lfloor n/3 \rfloor} = [a_i^{\lfloor n/3\rfloor},a_j] = 1$.  This completes the proof when $l = 2$.

We now assume that $l \ge 2$.  If $|G^l| \le p^{n/(l+1)}$, then the result holds.  Thus, we may assume that $|G^l| > p^{n/(l+1)}$, and so, $|G/G^l| < p^{nl/(l+1)}$.  We see that $G^l$ is generated by the set $\{[a,b] \mid a \in G^{l-1}, b \in G \}$.  As above, since $G^l$ is central (thus abelian) and a $p$-group, it suffices to show that $o([a,b]) \le p^{n/{(l+1)}}$ when $a \in G^{l-1}$ and $b \in G$.  By induction, we know since $G/G^l$ has nilpotence class $l-1$ that $o(a G^l)$ is less than or equal to $p^{(nl/(l+1))/l} = p^{n/(l+1)}$.  We then have that $[a,b]^{\lfloor n/(l+1) \rfloor} = [a^{\lfloor n/(l+1) \rfloor},b] = 1$.  This proves the desired result.
$\Box$\\ [1ex]

We now have what we need to prove the theorem.  Notice that when $l = 1$, $G$ is abelian and the inequality is the inequality that was proved in Theorem \ref{one}.  Thus, we can use the abelian case as the base case for our induction.  We do need two results that we have proved in \cite{pre1}: (1) If $N$ is a normal subgroup of $G$, then $\eta (G) \ge \eta (G/N)$. (2) If $G$ is a group, then $\eta (G) \ge \eta (Z(G))$.


\bigskip

{\bf Proof.} [Proof of Theorem \ref{nilp cl}]
We work by induction on $l$.  If $l = 1$, then this is Theorem \ref{one}.

We now suppose that $l \ge 2$.  Note that $G/G' = G/G^2$ is not cyclic, so the induction hypothesis is valid even if $l =2$.  If $|G/G^l| \ge p^{n(l-1)/l}$, then by the induction hypothesis we have 
$$\eta (\frac G{G^l}) \ge (p-1) \left(\frac {n(l-1)/l}{l-1} - 2\right) + p + 1 = (p-1)\left(\frac nl - 2\right) + p + 1.$$
As we noted before the start of the proof, we have that $\eta (G) \ge \eta (G/G^l)$, and so we have the conclusion in this case.  

Thus, we may assume that $|G/G^l| < p^{n(l-1)/l}$, and so, $|G^l| > p^{(n/l)} > p^{n/(l+1)}$.  By Lemma \ref{orders}, this implies that $G^l$ is not cyclic.  It follows that $|Z(G)| > p^{n/l}$ and is not cyclic.  By Theorem \ref {one}, we have $\eta (Z(G)) \ge (p-1)(n/l - 2) + p + 1$, and 
as we noted before the start of the proof, we have $\eta (G) \ge \eta (Z(G))$.  This yields the result when $n = l$ and proves the theorem.
$\Box$\\ [1ex]

Rewriting the inequality in Theorem \ref{nilp cl}, we obtain $\ln_p (|G|) \le ((\eta (G) - p -1)/(p-1) + 2)l$; so Theorem \ref{nilp cl} can be viewed as saying that that the size of $G$ is bounded in terms of a function of $\eta (G)$ and the nilpotence class of $G$.  I.e., if the order of $G$ is growing, then either $\eta$ or $l$ must be growing.  If one fixes $|G|$, then the nilpotence class being large forces $\eta$ to be small and vice versa.  In this way, it is like the uncertainty principle.



\noindent Mariagrazia Bianchi:
Dipartimento di Matematica F. Enriques,
Universit\`a degli Studi di Milano, via Saldini 50,
20133 Milano, Italy.\\
mariagrazia.bianchi@unimi.it\\[1ex]
Rachel D. Camina: Fitzwilliam College, Cambridge, CB3 0DG, UK.\\
rdc26@cam.ac.uk\\[1ex]
Mark L. Lewis:  Department of Mathematical Sciences, Kent State University, Kent, Ohio, 44242 USA.\\
lewis@math.kent.edu\\


\begin{thebibliography}{99}
	
	
	
\bibitem{pre1}
\rm M. Bianchi, R. D. Camina, Mark L. Lewis, E. Pacifici, Conjugacy classes of maximal cyclic subgroups, Preprint.
	
\bibitem{pre}	
\rm Yiftach Barnea, Rachel D. Camina, Mikhail Ershov, Mark L. Lewis, Preprint.

\bibitem{east}
\rm T. E. Easterfield, The orders of products and commutators in prime-power groups, Proc. London Math. Soc. (2), {\bf 36} (1933), 73-77.

	
	
	
	
	
	
\bibitem{rogerio}
\rm J. R. Rog{\'e}rio, A note on maximal coverings of groups, Comm. Algebra {\bf 42}(10) (2014), 4498-4508.
	
	
	
\bibitem{von}
\rm Timm Wilhelm von Puttkamer,  On the Finiteness of the Classifying Space for Virtually Cyclic Subgroups,
PhD thesis.
	
	
\bigskip	
	
	
	
\end{thebibliography}
\end{document}